\documentclass[11pt]{article}

\usepackage{latexsym,amssymb,amsmath,amsthm}
\usepackage{graphicx}
\usepackage{enumerate}
\numberwithin{equation}{section} 

\theoremstyle{definition}
\newtheorem{definition}{Definition}[section]

\theoremstyle{remark}
\newtheorem{remark}{Remark}[section]

\theoremstyle{plain}
\newtheorem{proposition}{Proposition}[section]
\newtheorem{corollary}{Corollary}[section] 
\newtheorem{lemma}{Lemma}[section]
\newtheorem{theorem}{Theorem}[section]

\begin{document}
  \title{Two-Parameter L\'evy Processes \\ along Decreasing Paths}
   \author{Shai Covo
   \footnote{A preprint version of this paper forms part of the author's PhD thesis, prepared at Bar Ilan University under the supervision of Prof. E. Merzbach. Work supported by the Doctoral Fellowship of Excellence,
             Bar Ilan University.}\\  
   Department of Mathematics\\Bar Ilan University\\52900 Ramat-Gan, Israel\\\small E-mail: green355@netvision.net.il} 
\date{\today}
\maketitle

\begin{abstract} 
Let $\{ X_{t_1,t_2} : t_1,t_2 \geq 0 \}$ be a two-parameter L\'evy process on $\mathbb{R}^d$.
We study basic properties of the one-parameter process $\{ X_{x(t),y(t)} :t \in T\}$ where $x$ and $y$ are, respectively, nondecreasing and nonincreasing nonnegative continuous functions on the interval $T$. We focus on and characterize the case where the process has stationary increments. 
\end{abstract}
  
  \textsl{keywords}: two-parameter L\'evy processes; decreasing paths; stationary increments; functional equation; Brownian sheet; Brownian bridge.

\section{Preliminaries}

In this paper, we study a class of one-parameter processes obtained from two-parameter L\'evy processes (TPLPs, defined below) by restricting the latter to paths (parameterized curves) in $\mathbb{R}_ + ^2 $ as indicated in the abstract. This topic---which seems to have received almost no attention in the literature---is vast; the present work provides a solid basis for further research. 

We first collect from \cite{dalang2,ehm,lagaize,vares} some basic material on TPLPs. These processes are often called {\it L\'evy sheets}, in analogy with the Brownian case. Their extension to $n>2$ parameters is straightforward (see e.g. \cite{ehm}).

TPLPs are indexed by ${\mathbb R}_ + ^2  = [0,\infty)^2$. A typical parameter (``time point'') $t \in {\mathbb R}_ + ^2$ is written as $t=(t_1,t_2)$. For times $s,t \in {\mathbb R}_ + ^2$ with $s_1 \leq t_1$ and $s_2 \leq t_2$, we write $(s,t]$ for the {\it rectangle} $(s_1,t_1] \times (s_2,t_2]$.
For a function $f:\mathbb{R}^{2}_{+} \rightarrow \mathbb{R}^d$, the {\it increment} $f((s,t])$ of $f$ over the rectangle $(s,t]$ is defined to be
\begin{equation} \nonumber 
f((s,t]) = f(t_1 ,t_2) - f(s_1 ,t_2) - f(t_1 ,s_2) + f(s_1 ,s_2).
\end{equation}
(Note that $f(\emptyset)=0$.) For a finite collection of disjoint rectangles $B_1,\ldots,B_n$, we can then define
\begin{equation} \nonumber
f\bigg(\mathop  \bigcup \limits_{i = 1}^n B_i \bigg) = \sum\limits_{i = 1}^n {f(B_i )}. 
\end{equation}
For a point $t \in {\mathbb R}_ + ^2$, we determine the following quadrants (as in \cite{lagaize}):
\begin{equation} \nonumber
\begin{split}
Q_1 (t) & = \{ s \in \mathbb{R}_ + ^2 :s_1  \ge t_1 ,\,s_2  \ge t_2 \} ,\;\;Q_2 (t) = \{ s \in {\mathbb R}_ + ^2 :s_1  < t_1 ,\,s_2  \ge t_2 \} ,\\
Q_3 (t) & = \{ s \in {\mathbb R}_ + ^2 :s_1  < t_1 ,\,s_2  < t_2 \} ,\;\; Q_4 (t) = \{ s \in {\mathbb R}_ + ^2 :s_1  \ge t_1 ,\,s_2  < t_2 \} .
\end{split}
\end{equation}
The function $f$ is said to have {\it quadrantal limits} if for each point $t \in {\mathbb R}_ + ^2 $ and $i=1,\ldots,4$, the four limits
\begin{equation} \nonumber
f (t^{(i)}): = \mathop {\lim }\limits_{\scriptstyle s \rightarrow t \atop 
  \scriptstyle s \in Q_i (t)} f(s)
\end{equation}
exist whenever $Q_i (t) \neq \emptyset$. It is said to be {\it right continuous} if $f(t)=f(t^{(1)})$ for all $t$. 
A function $f$ right continuous with quadrantal limits is continuous except on at most countably many horizontal and vertical lines (see e.g. \cite[p. 163]{vares}).
For such $f$, the {\it jump} of $f$ at $t \in (0,\infty)^2$ is defined to be 
\begin{equation} \nonumber
\begin{split}
\Delta f(t) & = f(t^{(1)}) - f(t^{(2)} ) - f(t^{(4)} ) + f(t^{(3)} ) \\
            & = f(t_1 ,t_2 ) - f(t_1^ -  ,t_2 ) - f(t_1 ,t_2^ -  ) + f(t_1^ -  ,t_2^ -  ).
\end{split}
\end{equation}
If $f$ vanishes on the axes, we set $\Delta f(t)=0$ for any $t \in \mathbb{R}_+^2$ with $t_1t_2=0$.

The above definitions are applied to a two-parameter process $X$ upon identifying its sample paths
$X_t (\omega ):\mathbb{R}^{2}_{+} \rightarrow \mathbb{R}^d$ with $f$.

\begin{definition} \label{dfn1}
A stochastic process $X = \{ X_t: t \in {\mathbb R}_ + ^2 \}$ taking values in $\mathbb{R}^d$ is a {\it two-parameter L\'evy process} if the following conditions are satisfied. 
\begin{enumerate}[{\rm (i)}]

\item For any choice of $n \geq 2$ and disjoint rectangles $B_1,\ldots,B_n$, the random variables $X(B_1),\ldots,X(B_n)$ are independent.

\item If $B$ is a rectangle and $p \in \mathbb{R}^{2}_{+}$, then $X(B)$ and $X(B+p)$ are identically distributed, where $B+p=\{t+p:t \in B\}$.

\item $X$ vanishes on the axes a.s. (almost surely).

\item $X$ is continuous in probability.

\item The sample paths of $X$ are a.s. right continuous with quadrantal limits.

\end{enumerate}
\end{definition} 

(Any process satisfying (i)-(iv) may be called a {\it TPLP in law}.)

The correspondence between infinitely divisible (ID) distributions and TPLPs can be stated as follows (see e.g. \cite[Theorem 1.1]{vares}). Let $\mu$ be an ID distribution on $\mathbb{R}^d$ and $\varphi$ its characteristic function. Then, there exists a TPLP $X$ such that $\varphi_{t_1,t_2} = \varphi^{t_1 t_2}$, for all $t \in \mathbb{R}^2_+$, where $\varphi_{t_1,t_2}$ is the characteristic function of $X_{t_1,t_2}$. Conversely, if $X$ is a TPLP, then the characteristic function of $X_{t_1,t_2}$ is of the above form. The law of a TPLP is determined by its one-dimensional distribution at time $t=(1,1)$.

Letting $\langle \cdot,\cdot \rangle$ denote the scalar product in $\mathbb{R}^d$,
the characteristic function of the increment $X(B)$ of a TPLP $X$ over a rectangle $B$ is thus given by   
\begin{equation} \label{eq:formula}
{\rm E}[{\rm e}^{{\rm i} \langle z,X(B)  \rangle } ] = {\rm e}^{m(B) \psi (z)} ,\;\; z \in \mathbb{R}^d,
\end{equation}
where here and in the sequel $m$ denotes Lebesgue measure on $\mathbb{R}^{2}_{+}$, and, by the L\'evy--Khintchine formula, $\psi:\mathbb{R}^d \rightarrow \mathbb{C}$ admits a unique representation 
\begin{equation} \label{eq:exponent_form}
\psi (z) =  {\rm i} \langle \gamma ,z \rangle  - \frac{1}{2} \langle z,Az \rangle + \int\nolimits_{\mathbb{R}^d } {({\rm e}^{{\rm i} \langle z,x \rangle }  - 1}  - {\rm i} \langle z,x \rangle \mathbf{1}_{\{ |x| \le 1\} }(x) )\nu ({\rm d}x),
\end{equation}
where $A$ is a symmetric nonnegative definite $d\times d$ matrix, $\nu$ is a measure on $\mathbb{R}^d$ satisfying $\nu (\{ 0\} ) = 0$ and 
$\int_{\mathbb{R}^d } {(|x|^2  \wedge 1)\nu ({\rm d}x)}  < \infty$, and $\gamma \in \mathbb{R}^d$. 
As in the one-parameter case, $\psi$ and $\nu$ are called the ({\it characteristic}) {\it exponent} and the {\it L\'evy measure} of $X$, respectively; if finite, $\gamma _0  = \gamma  - \int_{\{ |x| \le 1\} } {x \nu ({\rm d}x)}$ is called the {\it drift} of $X$; and, when $A=0$, $X$ is called {\it purely non-Gaussian}. 

\begin{remark} \label{rm:useful}
It will be useful to note that if $\varphi(z)={\rm e}^{f(z)}$ is the characteristic function of some ID random variable on $\mathbb{R}^d$, then $f$ corresponds to a characteristic exponent if and only if $f$ is continuous and $f(0)=0$; see Lemmas 7.5, 7.6 and the sentence after (8.9) in \cite{sato}.
We will also need the simple fact that if the one-dimensional distribution at time $t=(1,1)$ of a sequence $(X_n)$ of TPLPs on $\mathbb{R}^d$ converges to that of a TPLP $X$, then the finite-dimensional distributions (FDDs) of $X_n$ converge to those of $X$.  
\end{remark}

Like in the one-parameter case, a TPLP $X$ on $\mathbb{R}^d$ can be decomposed as 
\begin{equation} \nonumber
X_{t_1 ,t_2 }  =  t_1t_2 \gamma' + X_{t_1 ,t_2 }^0  + X_{t_1 ,t_2 }^1  + X_{t_1 ,t_2 }^2,\;\; t_1,t_2\geq 0,
\end{equation} 
where $X^i$ are independent TPLPs, such that $X^0$ is a continuous centered Gaussian process, $X^2$ is a compound Poisson process (CPP) with jumps (if any) of absolute value larger than some fixed $a>0$, $X^1$ has mean zero and jumps (if any) not exceeding $a$ in absolute value, and $\gamma' = \gamma'(a) \in \mathbb{R}^d$. 
More specifically, $X^1$ is obtained as the almost sure, uniform-on-compacts limit of a sequence of CPPs  `compensated' by having their means subtracted (a {\it compensated sum of jumps}). The L\'evy--Khintchine representation \eqref{eq:exponent_form} corresponds in an obvious way to the above decomposition with $a=1$. 

\begin{figure}
	\centering
	\includegraphics[height=50mm]{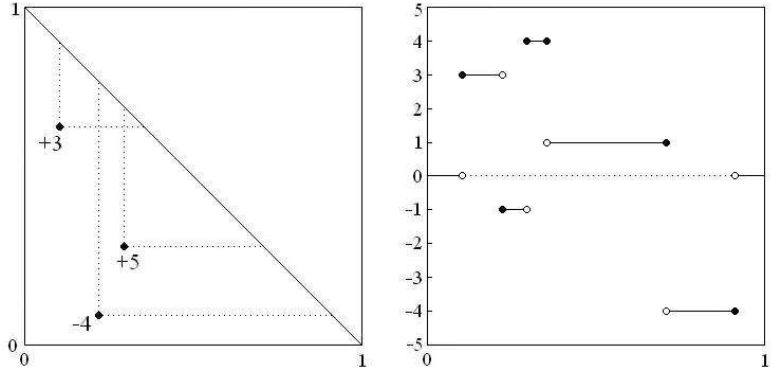}
	\caption {Jump locations and sizes of a two-parameter CPP $X$ in $\{ t \in {\mathbb R}_ + ^2 :t_1  + t_2  \le 1\}$, and the associated process $\{ X_{t,1 - t} :0 \le t \le 1\}$.}
	\label{fig1}
\end{figure}

To better understand the sample paths of TPLPs, one has to consider the set of their discontinuities. We note here that each discontinuity of $X$ propagates both horizontally and vertically from a jump of $X$. For specific details, we refer to \cite[Sect. 2.4]{dalang2}. See also Figure \ref{fig1} of the present paper.

We call an $n$-parameter (for $n \geq 1$ in general) L\'evy process $X$ {\it deterministic} (respectively, {\it zero}) if $X_t$ has a $\delta$ (respectively, $\delta_0$)-distribution for some (equivalently, any) point $t$ with positive coordinates. {\it Nondeterministic} and {\it nonzero} L\'evy processes are defined accordingly.
The process $X$ is {\it symmetric} if $X$ and $-X$ are identical in law, which amounts to say that $X_t \stackrel{{\rm d}}{=} -X_t$ for some (equivalently, any) point $t$ with positive coordinates. Here and in the sequel, $\stackrel{{\rm d}}{=}$ denotes equality in distribution of random variables/vectors.

Having introduced standard terminology, we now formally define the basic notions of the topic at hand.
Henceforth, we let $T  \subseteq \mathbb{R}$ denote some interval of the real line, and $T^\circ$ its interior. 

\begin{definition} \label{dfn2}
A {\it decreasing path} in $\mathbb{R}^{2}_{+}$ is a parameterized curve $\alpha(t)=(x(t),y(t))_{t \in T}$ where $x$ and $y$ are, respectively, nondecreasing and nonincreasing (continuous) functions on $T$, both strictly positive on $T^\circ$ and at least one is not identically constant.
\end{definition}
 
\begin{definition} \label{dfn3}
Given a TPLP $X$ on $\mathbb{R}^d$ and a decreasing path $\alpha(t)=(x(t),y(t))_{t \in T}$, we refer to the one-parameter process $\{ X_{x(t),y(t)} :t \in T\}$ as the {\it TPLP} $X$ {\it along the path} $\alpha$, and denote it by $X^\alpha$.
\end{definition} 

Figure \ref{fig1} illustrates a hypothetical sample path of a two-parameter CPP $X$ along the path $\alpha(t)=(t,1-t)_{t \in [0,1]}$. (Note that $X^\alpha$ is not right continuous with left limits.) As we have indicated above, each jump of $X$ inside the triangle $\{ t \in {\mathbb R}_ + ^2 :t_1  + t_2  \le 1\}$ gives rise to a discontinuity which propagates along the vertical and horizontal half-lines which emanate from the jump location. This results in two cancelling jumps in the sample path of $X^\alpha$. 

The remainder of this paper is organized as follows. In Sect. 2 we study general basic properties of TPLPs along decreasing paths, thus getting an impression of how rich and complex they are. In particular, we characterize their FDDs (Theorem \ref{thm1}). In Sect. 3 we consider separately the special, tractable Brownian sheet case. In Sect. 4 we state and prove the main result of this paper, Theorem \ref{thm2}, which fully characterizes the case where TPLPs have stationary increments along decreasing paths. The key to the proof is Lemma \ref{lem2}, stating the solution of some functional equation. In Sect. 5 we discuss the interesting classes of processes from Theorem \ref{thm2}.

\section{General basic properties} 

As we will see in the course of this paper, TPLPs along decreasing paths form a rich and intriguing class of one-parameter processes. We first show that these processes typically have dependent increments. Given a decreasing path $\alpha(t)=(x(t),y(t))_{t \in T}$, an integer $n \geq 1$, and times $t_1 , \ldots, t_n  \in T$ with $t_1<\cdots<t_n$, we introduce the following collection of {\it disjoint} rectangles:
\begin{equation} \nonumber
\{ B_{ij}  = (x(t_{i - 1} ),x(t_i )] \times (y(t_{i + j} ),y(t_{i + j - 1} )]:i = 1,\ldots,n,\;j = 1,\ldots,n - i + 1\}
\end{equation}
(some of them may be empty), with the convention, to be used also in Theorem \ref{thm1} below, that $x(t_0 ) = y(t_{n + 1} ) = 0$.
Letting $n=3$, we have
\begin{equation} \nonumber
\begin{split}
X_{t_3 }^\alpha   - X_{t_2 }^\alpha   & = [X(B_{31} ) - X(B_{12} )] - X(B_{21} ),\\
X_{t_2 }^\alpha   - X_{t_1 }^\alpha   & = [X(B_{22} ) - X(B_{11} )] + X(B_{21} ),
\end{split}
\end{equation}
for any TPLP $X$ on $\mathbb{R}^d$.
It follows that $X_{t_3 }^\alpha   - X_{t_2 }^\alpha$ and $X_{t_2 }^\alpha   - X_{t_1 }^\alpha$ are dependent if $m(B_{21}):=(x(t_2)-x(t_1))(y(t_2)-y(t_3))>0$ and $X$ is nondeterministic. We thus conclude that 
there are only three families of decreasing paths along which a {\it nondeterministic} TPLP has independent increments: 
\begin{itemize}
\item The horizontal paths ($y=$ const).
\item The decreasing vertical paths ($x=$ const).
\item The decreasing `first vertical, then horizontal' paths.
\end{itemize}
Moreover, the independent increments property corresponding to the last two cases means that
for any $n \geq 2$ and $t_1,\ldots,t_n \in T$ with $t_1<\cdots<t_n$, $X_{t_2 }^\alpha - X_{t_1 }^\alpha,\ldots,X_{t_n }^\alpha   - X_{t_{n - 1} }^\alpha $ are independent,
and is thus substantially weaker than the one that requires these random variables to be independent of $X_{t_1 }^\alpha$ as well. At this point, it is interesting to observe that if $\alpha(t)=(x(t),y(t))_{t > 0}$ is any decreasing path such that $x(t)y(t)=t$, then the one-dimensional distributions of $X^\alpha$ (for $t>0$) are identical to those of a one-parameter L\'evy process with the same exponent.

We proceed to characterize the individual increments of TPLPs along decreasing paths.
Given a decreasing path $\alpha$ and points $s,t \in T$ with $s<t$, we define the `upper' and `lower' rectangles $B_{s,t}^{\rm u}$ and $B_{s,t}^{\rm l}$ by
\begin{equation} \nonumber
B_{s,t}^{\rm u} = (0,x(s)] \times (y(t),y(s)]\;\;{\rm and} \;\;B_{s,t}^{\rm l} = (x(s),x(t)] \times (0,y(t)].
\end{equation}
Since, for any TPLP $X$ on $\mathbb{R}^d$,
\begin{equation} \nonumber
X_t^\alpha   - X_s^\alpha  = X(B_{s,t}^{\rm l} ) - X(B_{s,t}^{\rm u} )
\end{equation}
and $B_{s,t}^{\rm l}  \cap B_{s,t}^{\rm u}  = \emptyset$, it follows by independence and \eqref{eq:formula} that
\begin{equation} \label{eq:dif}
{\rm E}[{\rm e}^{{\rm i} \langle z,X_t^\alpha  - X_s^\alpha  \rangle }] = \exp[{m (B_{s,t}^{\rm l} )\psi (z) + m (B_{s,t}^{\rm u} )\psi ( - z)}],\;\; z \in \mathbb{R}^d.
\end{equation}
If $X$ is symmetric then $\psi(z)=\psi(-z)$, and so
\begin{equation} \label{eq:difsym}
\begin{split}
{\rm E}[{\rm e}^{{\rm i} \langle z,X_t^\alpha  - X_s^\alpha  \rangle } ] & = \exp \big[\big(m(B_{s,t}^{\rm l} ) + m(B_{s,t}^{\rm u})\big)\psi (z)] \\
& = \exp \big[\big(x(s)y(s) + x(t)y(t) - 2x(s)y(t)\big)\psi (z)\big].
\end{split}
\end{equation}
Equations \eqref{eq:dif}-\eqref{eq:difsym} are the starting point of the proof of our Theorem \ref{thm2}. 


We now characterize the FDDs (law) of TPLPs along decreasing paths.

\begin{theorem} \label{thm1}
Let $X$ be a TPLP on $\mathbb{R}^d$ with exponent $\psi$, and 
$\alpha(t)  = (x(t),y(t))_{t \in T}$ a decreasing path.
Then for any $n \geq 1$ and $t_1 , \ldots, t_n  \in T$ with $t_1<\cdots<t_n$, the characteristic function $\phi_{t_1,\ldots,t_n}$ of the $\mathbb{R}^{nd}$-valued random variable $(X_{t_k }^\alpha  )_{k = 1,\ldots,n}$ is given, for $z = (z_k )_{k = 1,\ldots,n}$, $z_k  \in \mathbb{R}^d$, by: 
\begin{multline} \label{eq:characteristic}
\phi _{t_1,\ldots,t_n} (z) \\
= \exp \Bigg[\sum\limits_{i = 1}^n {\sum\limits_{j = 1}^{n - i + 1} {(x(t_i ) - x(t_{i - 1} ))(y(t_{i + j - 1} ) - y(t_{i + j} ))\psi \Bigg(\sum\limits_{k = i}^{i + j - 1} {z_k } \Bigg)\Bigg]} }.
\end{multline}
\end{theorem} 

\proof
By definition,
\begin{equation} \nonumber
\phi _{t_1,\ldots,t_n } (z) = {\rm E}[{\rm e}^{{\rm i} \langle (z_k),(X_{t_k }^\alpha) \rangle } ] = 
{\rm E}\bigg[\exp \bigg({\rm i}\sum\limits_{k = 1}^n {\langle z_k ,X_{t_k }^\alpha  \rangle } \bigg)\bigg].
\end{equation}
With $B_{ij}$ as defined above, $X_{t_k }^\alpha$ can be decomposed as 
\begin{equation} \nonumber
X_{t_k }^\alpha   = \sum\limits_{i = 1}^k {\sum\limits_{j = k - i + 1}^{n - i + 1} {X(B_{ij} )} }. 
\end{equation}
Noting that the term $X(B_{ij} )$, for fixed $i$ and $j$, appears in this decomposition for $k = i,\ldots,i + j - 1$ only, we get by virtue of independence
\begin{equation} \nonumber
{\rm E}\bigg[\exp \bigg({\rm i}\sum\limits_{k = 1}^n {\langle z_k ,X_{t_k }^\alpha  \rangle } \bigg)\bigg] = 
\prod\limits_{i = 1}^n {\prod\limits_{j = 1}^{n - i + 1} {{\rm E}\bigg[\exp \bigg({\rm i}\bigg \langle \sum\limits_{k = i}^{i + j - 1} {z_k } ,X(B_{ij} )\bigg \rangle } \bigg)\bigg]}.
\end{equation}
The proof then follows from \eqref{eq:formula} and the value of $m (B_{ij} )$. \hfill $\Box$ \\  
 
\begin{corollary} \label{cor1}
Let $X$ be a nondeterministic TPLP on $\mathbb{R}^d$, and 
$\alpha _i (t) = (x_i (t),y_i (t))_{t \in T}$, $i = 1,2$, decreasing paths. Then $X^{\alpha _1 } \stackrel{{\rm law}}{=} X^{\alpha _2 }$ if and only if 
$x_2 (t) = p x_1 (t)$ and $y_2 (t) = p^{-1}y_1 (t)$
for all $t \in T$, for some $p > 0$.
\end{corollary}

\proof
The ``if'' part is clear from Theorem \ref{thm1}. 
For the ``only if'' part, we assume that $(X_{t_1 }^{\alpha _1 } ,X_{t_2 }^{\alpha _1 } ) \stackrel{{\rm d}}{=} (X_{t_1 }^{\alpha _2 } ,X_{t_2 }^{\alpha _2 } )$, for all $t_1,t_2 \in T$ with $t_1<t_2$, and in accordance with \eqref{eq:characteristic} express the identical characteristic functions as
\begin{equation} \label{eq:characteristics}
\phi _{i;t_1 ,t_2}{(z_1,z_2)}  = \exp \big[a_{i} \psi (z_1 ) + b_{i} \psi (z_1  + z_2 ) + c_{i} \psi (z_2 )\big],\;\;i=1,2,
\end{equation}
respectively. Noting that $\psi(0)=0$, it follows straightforwardly that $a_1  + b_1  = a_2  + b_2$ and $b_1  + c_1  = b_2  + c_2$. Hence, $a_1-a_2 = c_1-c_2$.
Further, setting $z_1= z$ and $z_2=- z$ in \eqref{eq:characteristics} leads to $(a_1  - a_2 )\psi ( z) = (c_1  - c_2 )( - \psi ( -  z))$. By the assumption on $X$, we have $\psi ( z) \ne  - \psi ( - z)$ for some $ z \in \mathbb{R}^d$; see the proof of Lemma \ref{lem3} below.  Combining it all, we conclude that $b_1  = b_2$, and explicitly
$x_1 (t_1 )y_1 (t_2 ) = x_2 (t_1 )y_2 (t_2 )$. The corollary is thus established. \hfill $\Box$ \\ 

The following lemma, well known for Poisson random variables, is of general interest. 

\begin{lemma} \label{lem1}
Let $\xi_i$, $i=1,2$, be independent, ID, integrable random variables on $\mathbb{R}$ with characteristic functions ${\rm E}[{\rm e}^{{\rm i}z\xi _i } ] = {\rm e}^{c_i \psi (z)}$, $c_i>0$ fixed. Then,
\begin{equation} \nonumber
{\rm E}[\xi _1 |\xi _1  + \xi _2 ] = \frac{{c_1 }}{{c_1  + c_2 }}(\xi _1  + \xi _2 ).
\end{equation}
Put another way, given an integrable L\'evy process $Y$ on $\mathbb{R}$, it holds that
\begin{equation} \nonumber
{\rm E}[Y_s | Y_t] = \frac{s}{t}Y_t,\;\;0<s<t.
\end{equation}
\end{lemma}

\proof We prove the second formulation. Assume first that $s/t=m/n$, with $m,n \in \mathbb{N}$. From $\sum\nolimits_{i = 1}^n {{\rm E}[Y_{it/n}  - Y_{(i - 1)t/n} |Y_t ]}  = Y_t$ we deduce that ${\rm E}[Y_{t/n}|Y_t]=Y_t / n$, and, in turn, ${\rm E}[Y_s |Y_t ] = (m/n)Y_t  = (s/t)Y_t$. If $s/t$ is irrational, let $(s_j)$ be a sequence such that $s_j  \uparrow s$ with $s_j/t$ being rational. By an elementary property of L\'evy processes, $Y_{s_j } \stackrel{{\rm a.s.}}{\rightarrow} Y_s $. Define $Y_s^*  = \sup _{u \in [0,s]} |Y_u |$; thus $|Y_{s_j}|\leq Y_s^*$ $\forall j$. Since, by assumption, ${\rm E}[|Y_s|]<\infty$, we conclude \cite[Theorem 25.18]{sato} that also ${\rm E}[Y_s^*]<\infty$. Hence, by the dominated convergence theorem for conditional expectations, ${\rm E}[Y_{s_j } |Y_t ] \stackrel{{\rm a.s.}}{\rightarrow} {\rm E}[Y_s |Y_t ]$. Since $s_j/t$ is rational, ${\rm E}[Y_{s_j } |Y_t ]=(s_j/t)Y_t$. Thus, ${\rm E}[Y_s |Y_t ] = (s/t)Y_t$. \hfill $\Box$ \\

\begin{corollary} \label{cor2}
Let $X$ be an integrable TPLP on $\mathbb{R}$ and 
$\alpha(t)  = (x(t),y(t))_{t \in T}$ a decreasing path. Then, for every $s,t \in T$ with $s<t$,
\begin{equation} \label{eq:general_expectation}
{\rm E}[X_t^\alpha  |X_s^\alpha  ] = \frac{{y(t)}}{{y(s)}}X_s^\alpha   + (x(t) - x(s))y(t){\rm E}[X_{1,1} ].
\end{equation}
In particular, if $X$ is moreover nondeterministic and if there exist $s,t \in T$, $s<t$, with $x(s)<x(t)$ and $y(t)<y(s)$ (i.e., the path image is not vertical or horizontal), then $X^\alpha$ is neither a submartingale nor a supermartingale.
\end{corollary}  

\proof Equation \eqref{eq:general_expectation} follows straightforwardly by writing $X_t^\alpha$ and $X_s^\alpha$ as
\begin{equation} \nonumber
\begin{split}
X_t^\alpha   & = X((0,x(s)] \times (0,y(t)]) + X(B_{s,t}^{\rm l} ), \\
X_s^\alpha   & = X((0,x(s)] \times (0,y(t)]) + X(B_{s,t}^{\rm u} ),
\end{split}
\end{equation}
and applying Lemma \ref{lem1}. The second part of the corollary follows easily. \hfill $\Box$ \\ 

\begin{remark} \label{rm:intuitive}
Despite its importance, we will not address the question of markovity of TPLPs along decreasing paths. However, it is intuitively clear that such processes are generally non-Markovian (consider, e.g., Figure \ref{fig1}).
\end{remark}

\section{The Brownian case}

A two-parameter standard Brownian sheet on $\mathbb{R}^d$ is a TPLP with exponent 
\begin{equation} \nonumber 
\psi (z) =  - \frac{1}{2}|z|^2 ,\;\;z \in \mathbb{R}^d.
\end{equation}
We shall denote it by $W$. Equivalently, $W$ is a continuous centered Gaussian process indexed by $\mathbb{R}_+^2$ and taking values in $\mathbb{R}^d$, with covariance structure given by the following: for all $s,t \in \mathbb{R}_+^2$ and
$1 \leq i,j\leq d$,
\begin{equation} \nonumber
{\rm E}[W_{s_1 ,s_2 }^i W_{t_1 ,t_2 }^j ] = \begin{cases}
(s_1  \wedge t_1 )(s_2  \wedge t_2 ),& {\rm if}\; i=j; \\
0,& {\rm if}\; i \neq j,
\end{cases}
\end{equation}
where the superscripts refer to vector components. Put another way, $W = (W^i )_{i=1,\ldots,d}$
where $W^i$ are independent standard real-valued Brownian sheets. 

Given a decreasing path $\alpha (t) = (x(t),y(t))_{t \in T}$, the process $W^\alpha := \{ W_{x(t),y(t)} :t \in T\}$ is thus an $\mathbb{R}^d$-valued continuous centered Gaussian process with covariance structure given for $s,t \in T$ by
\begin{equation} \label{eq:covariance_decreasing}
{\rm E}[(W_s^\alpha)^i (W_t^\alpha)^j ] = \begin{cases}
x(s \wedge t)y(s \vee t),& {\rm if}\; i=j; \\
0,& {\rm if}\; i \neq j.
\end{cases}
\end{equation}
Since the law of a centered Gaussian process is determined by its covariance structure, Corollary \ref{cor1} for the case $X=W$ follows immediately from \eqref{eq:covariance_decreasing}.

Rather than using \eqref{eq:characteristic}, we express the characteristic function $\phi_{t_1,\ldots,t_n}$ of the $\mathbb{R}^{nd}$-valued random variable $(W_{t_i }^\alpha  )_{i = 1,\ldots,n}$, in a more elegant form, using \eqref{eq:covariance_decreasing} and the standard formula for the characteristic function of a multivariate Gaussian random variable. Specifically, we have $\phi _{t_1 , \ldots ,t_n } (z) = \exp[{ - {\textstyle{1 \over 2}}\langle (z_i ),A(z_i )\rangle }]$, where $z = (z_i )_{i = 1, \ldots ,n}$, $z_i \in \mathbb{R}^{d}$, and $A$ is the $nd \times nd$ matrix given by $A = (C_{ij} )_{1 \le i,j \le n}$ with $C_{ij}=x(t_i \wedge t_j)y(t_i \vee t_j)I_d$, $I_d$ being the identity matrix of order $d$, which leads to 
\begin{multline} \label{eq:joint_normal}
\phi _{t_1,\ldots,t_n} (z) \\
= \exp \bigg[ - \frac{1}{2}\bigg(\sum\limits_{i = 1}^n {x(t_i )y(t_i )|z_i |^2 }  + 2\sum\limits_{i = 1}^{n - 1} {\sum\limits_{j = i + 1}^n {x(t_i )y(t_j )\langle z_i ,z_j \rangle } } \bigg)\bigg]
\end{multline}
(for any $n\geq1$ and $t_1 , \ldots, t_n  \in T$ with $t_1<\cdots<t_n$).


The following proposition reveals an appealing feature of Brownian sheets along decreasing paths. It is verified by comparing covariance structures. Henceforth, `Brownian motion' is abbreviated `BM'.

\begin{proposition} \label{prop2}
Let $\alpha(t)=(x(t),y(t))_{t \in T}$ be a decreasing path and $B$ a standard BM on $\mathbb{R}^d$ (i.e., $B = (B^i)_{i=1,\ldots,d}$ where $B^i$ are independent standard real-valued BMs). Then the following identities in law hold:
\begin{eqnarray} 
W^\alpha & \stackrel{{\rm law}}{=} & \{ y(t)B_{x(t)/y(t)}:{t \in T}\} \label{eq:idl1} \\ 
& \stackrel{{\rm law}}{=} & \{ x(t)B_{{y(t)/x(t)}}:{t \in T} \} \label{eq:idl2} \\
& \stackrel{{\rm law}}{=} & \{ (x(t) + y(t))B_{x(t)/(x(t) + y(t))} - x(t)B_1 :t \in T\} \label{eq:idl3} \\
& \stackrel{{\rm law}}{=} & \{ (x(t) + y(t))B_{y(t)/(x(t) + y(t))}  - y(t)B_1 :t \in T\}, \label{eq:idl4} 
\end{eqnarray}
where the right-hand processes are defined to be zero whenever $x(t)y(t)=0$.
\end{proposition}

This proposition provides a useful tool to simulate and analyze Brownian sheets along decreasing paths. (Note that the time indices of the BM in \eqref{eq:idl3}-\eqref{eq:idl4} are restricted to $[0,1]$.) The equality in law  \eqref{eq:idl1} also shows that the process $W^\alpha$ is a Gauss--Markov process (i.e., both Gaussian and Markovian). In the well-behaved case, the covariance function of an arbitrary real-valued Gauss--Markov process having no fixed values (singularities) in $T^\circ$ satisfies $c(s,t)=h_1{(s)}h_2{(t)}$, $s\leq t$, $s,t \in T^\circ$, where $r(t)=h_1{(t)}/h_2{(t)}$ is a positive monotonically increasing function on $T^\circ$ \cite[p. 455]{di-nardo}. A centered Gaussian process with such covariance can thus be represented (on $T^\circ$) as $h_2{(t)}B_{r(t)}$, with $B$ a standard BM. We conclude that Brownian sheets along decreasing paths constitute a fundamental class of (centered) Gauss--Markov processes. We now illustrate the importance of relation \eqref{eq:idl1}.

Suppose that $W$ is real-valued and $\alpha$ is a decreasing path such that $r(t):=x(t)/y(t)$ is strictly increasing on $T^\circ$. The transition density function $f(\cdot,t|z,s)$ ($s,t \in T$, $s<t$) of $W^\alpha$ can be found using \eqref{eq:idl1}, or immediately from the general result \cite[equations (2.5)]{di-nardo} for Gauss--Markov processes, to be the normal density with mean $\mu$ and variance $\sigma^2$ given by:
\begin{equation} \nonumber
\begin{split}
\mu & =  {\rm E}[W_t^\alpha  |W_s^\alpha   = z]  = \frac{{y(t)}}{{y(s)}}z, \\
\sigma^2  & =  {\rm Var}[W_t^\alpha  |W_s^\alpha   = z] = y(t)\bigg[x(t) - \frac{{y(t)}}{{y(s)}}x(s)\bigg]
\end{split}
\end{equation}
(recall \eqref{eq:general_expectation}), with the obvious interpretation in case $y(t)=0$.

With $\alpha$ and $r$ as in the last paragraph, let $P_0^z (r(s),r(t))$ denote the probability that a standard real-valued  BM $B$ has at least one zero in the time interval $(r(s),r(t))$, $s,t \in T^\circ$, given that $B_{r(s)}=z \neq 0$. Similarly, let $\tilde P_0^z (s,t)$ denote the probability that $W^\alpha$ has at least one zero in the time interval $(s,t)$ given that $W_s^\alpha = z\,(\neq 0)$. It follows readily from \eqref{eq:idl1} that $\tilde P_0^z (s,t) = P_0^{z/y(s)} (r(s),r(t))$. Hence, using an elementary formula for BM, 
\begin{equation} \nonumber
\tilde P_0^z (s,t) = \frac{{|z/y(s)|}}{{\sqrt {2\pi } }}\int_0^{r(t) - r(s)} {u^{ - 3/2} \exp \bigg[ - \frac{{(z/y(s))^2 }}{{2u}}\bigg] \,{\rm d}u}. 
\end{equation}
Similarly, using another elementary formula for BM, if we let $\tilde P_0 (s,t)$ denote the probability that $W^\alpha$ has at least one zero in the time interval $(s,t)$, then $\tilde P_0 (s,t) = 2\pi^{-1} \arccos \sqrt {r(s)/r(t)}\; (= 1- 2\pi^{-1}\arcsin \sqrt {r(s)/r(t)})$.


In light of \eqref{eq:idl1}, the following remark is in order.

\begin{remark} \label{rm:idl1}
Let $Y$ be a one-parameter L\'evy process on $\mathbb{R}^d$ and $x(\tau),\,y(\tau)$ any positive numbers with $y(\tau)\neq 1$. Then $Y_{x(\tau)y(\tau)} \stackrel{{\rm d}}{=} y(\tau)Y_{x(\tau)/y(\tau)}$ if and only if $Y$ is Gaussian with mean $0$. Since this is trivial if $Y$ is deterministic, we will assume it is not.

\proof Fix $x(\tau)$ and $y(\tau)$ as above, and put $a=y^2{(\tau)}$, $b=y(\tau)$, and $t'=x(\tau)/y(\tau)$. Since $Y_{at'}=Y_{x(\tau)y(\tau)}$ and $bY_{t'}=y(\tau)Y_{x(\tau)/y(\tau)}$, it suffices to show that the L\'evy processes $\{ Y_{at} :t \ge 0\}$ and $\{ bY_t :t \ge 0\}$ are identical in law if and only if $Y$ is Gaussian with mean $0$. Noting that $b=a^{1/2}\neq 1$, this equality in law reads ``the process $Y$ is {\it semi-selfsimilar} with {\it exponent} $H=1/2$'' \cite[Definitions 13.4, 13.12]{sato}. In turn, $\alpha=1/H=2$ is the {\it index} of $Y$ as a {\it strictly semi-stable} process (see Definition 13.16 and moreover Proposition 13.5 in \cite{sato}). From \cite[Theorem 14.2]{sato} we know that $Y$ is strictly semi-stable with index $2$ if and only if it is Gaussian with mean $0$. \hfill $\Box$ \\
\end{remark}       

We now turn to consider briefly the two particularly interesting examples of Brownian sheets along decreasing paths---both indicated, in normalized form, in \cite[Sect. 2.5]{dalang1}---namely, the Brownian bridge and OU processes. As usual, $B$ denotes a standard BM; $\alpha$ stands for a decreasing path.

An $\mathbb{R}^d$-valued standard Brownian bridge on $[0,l]$ is merely defined as a $d$-dimensional vector whose components are independent real-valued standard Brownian bridges on $[0,l]$. Losing nothing essential by assuming that $d=1$, we let $\hat B^l$ denote a standard real-valued Brownian bridge on $[0,l]$, that is, a centered Gaussian process on $[0,l]$ with covariance ${\rm E}[\hat B_s^l \hat B_t^l ] = s(1 - t/l)$, $0 \leq s \leq t \leq l$. Comparing covariances, we see that $W^\alpha \stackrel{{\rm law}}{=} \hat B^l$ if and only if $\alpha (t)  = (pt,p^{ - 1} (1 - t/l))_{t \in [0,l]}$ for some $p>0$. Setting $l=1$ and $p=1$, we thus recover from  \eqref{eq:idl3} the usual identity in law $\hat B^1 \stackrel{{\rm law}}{=}  \{B_t - tB_1:t \in [0,1]\}$ and from  \eqref{eq:idl1} the well-known one $\hat B^1 \stackrel{{\rm law}}{=} \{ (1-t)B_{t/(1-t)}:t \in [0,1]\}$ ($0B_\infty:=0$); 
from each of the counterpart equations \eqref{eq:idl4} and \eqref{eq:idl2} follows the elementary fact that $\hat B^1$ is invariant under the time reversal $t \mapsto 1-t$. The identity in law $\hat B^l  \stackrel{{\rm law}}{=} \{ W_{t,1 - t/l} :t \in [0,l]\}$ will be used in Sect. 5.2 to show that a Brownian bridge can be viewed as a difference of two limiting BMs.

We recall that a real-valued, centered, stationary Ornstein--Uhlenbeck process $V=\{ V_t:t\geq 0\}$ is a centered Gaussian process with covariance ${\rm E}[V_s V_t ] = r{\rm e}^{ - c|t - s|}$, $s,t\geq 0$, with $r$ and $c$ positive constants. Thus $W^\alpha \stackrel{{\rm law}}{=} V$ if and only if $\alpha (t) = (a{\rm e}^{ct} ,b {\rm e}^{ - ct} )_{t \ge 0}$ with $ab=r$. Equation \eqref{eq:idl1} leads to a common representation of $V$.
Let us set the variance parameter $r$ equal to $1/2$. 
Then, the process $V$ is representable as
\begin{equation} \label{eq:integral0}
V_t  = {\rm e}^{ - ct} V_0  + \int_0^t {{\rm e}^{ - c(t - s)} \,{\rm d}B_{cs} },
\end{equation}
where the standard BM $B$ is independent of $V_0$ (cf. e.g. \cite[Sect. 2]{valdivieso}).
If we take $B$ in \eqref{eq:integral0} to be a standard BM on $\mathbb{R}^d$, and accordingly $V_0$ a vector of $d$ i.i.d. ${\rm N}(0,1/2)$ variables (independent of $B$), we get an ordinary stationary OU process on $\mathbb{R}^d$.
In Sect. 5.3 we will consider the general case of stationary {\it processes of Ornstein--Uhlenbeck type}, where a L\'evy process $Z$ on $\mathbb{R}^d$ takes the place of $B$ in \eqref{eq:integral0}.

Finally, the following observation is in order.
Like the Brownian bridge and the OU process from above, a fractional BM with Hurst parameter $\neq 1/2$ is a continuous centered Gaussian process with stationary but dependent increments and is not a martingale (recall here the conclusion of Corollary \ref{cor2}). However, unlike its counterparts, it is not Markovian. Thus, it cannot
 be represented as $W$ along a decreasing path. Alternatively, this conclusion is an immediate consequence of Theorem \ref{thm2}, below.

\section{The main result}

The following lemma is the key to the proof of our main result, Theorem \ref{thm2}. 
It is an important result in its own right.

\begin{lemma} \label{lem2}
A decreasing path $(x(t),y(t))_{t \in T}$ satisfies the functional equation
\begin{equation} \label{eq:func_eq}
x(s)y(s) + x(t)y(t) - 2x(s)y(t) = \varphi (t - s),\;\; s,t \in T,  \; s<t
\end{equation}
for some (necessarily nonnegative) function $\varphi$ defined on the interval
$D = \{ t - s:s,t \in  T ,\; t > s\}$ if and only if one of the following holds:

\begin{enumerate}[{\rm (i)}]

\item $x(t) = a,\;y(t) = b - ct\;\; {\rm or}\;\; x(t) = b + ct,\;y(t) = a,\; t \in T$, for some $b \in \mathbb{R}$ and positive constants $a,c$. In both cases, $\varphi (u) = acu,\;u \in D$.

\item $x(t) = a + d(t - s^* )\mathbf{1}_{\{ t > s^* \} }(t),\;y(t) = b - c(t - s^* )\mathbf{1}_{\{ t \le s^* \} }(t),\; t \in T$, for some $s^* \in T^\circ$ and positive constants $a,b,c,d$ satisfying $ac=bd$. In this case, $\varphi (u) = acu,\;u \in D$.

\item $x(t) = a + bt,\;y(t) = c - dt,\; t \in T$, for some $a,c \in \mathbb{R}$ and positive constants $b,d$. In this case, $\varphi (u) = (ad + bc)u - bdu^2,\;u \in D$.

\item $x(t) = a{\rm e}^{ct},\;y(t) = b{\rm e}^{ - ct} ,\;t \in T$, for some positive constants $a,b,$ and $c$. In this case, $\varphi (u) = 2ab(1 - {\rm e}^{ - cu} ),\;u \in D$.
\end{enumerate}

\end{lemma} 

\proof
Sufficiency is easily verified. Necessity. Suppose that $(x(t),y(t))_{t \in T}$ is a decreasing path satisfying
\eqref{eq:func_eq}. Losing no generality, we assume that $T$ is {\it open}. In view of (i) above, it suffices to show that one of (ii)-(iv) must hold if 
both $x$ and $y$ belong to the class $\mathfrak{S}$ of strictly positive functions on $T$ which are not identically constant, which we henceforth assume.

The first key observation is that, by virtue of \eqref{eq:func_eq} and the almost everywhere differentiability of (the monotonic) $x$ and $y$ on $T$, $\varphi$ is twice differentiable on $D$. Then, isolating $x(t)$ as well as $y(s)$ in \eqref{eq:func_eq}, we have that 
\begin{equation} \nonumber
H: = \{ t \in T:\exists x''(t)\}  = \{ t \in T:\exists y''(t)\}. 
\end{equation}
From this and the following pair of equations,
\begin{equation} \nonumber
\begin{split}
(y(\tilde t) - y(t))x(s) & = \frac{1}{2}[\varphi (t - s) - \varphi (\tilde t - s) - x(t)y(t) + x(\tilde t)y(\tilde t)],\;\;\;s<t<\tilde t, \\
(x(\tilde s) - x(s))y(t) & = \frac{1}{2}[\varphi (t - s) - \varphi (t - \tilde s) - x(s)y(s) + x(\tilde s)y(\tilde s)],\;s<\tilde s<t,
\end{split}
\end{equation}
it follows easily that $s^* \in T\backslash H$ implies (ii) of the lemma.  

It remains to show that $H=T$ implies (iii) or (iv) of the lemma.
Suppose for a contradiction that $x'$ is not strictly positive on $T$. Then, since $x \in \mathfrak{S}$, there exists some $\tilde s \in T$ such that $x'(\tilde s)=0$ and $\{ s \in T:|s - \tilde s| < \delta ,\,x'(s) > 0\}$ is nonempty for any $\delta>0$. From \eqref{eq:func_eq} we have
\begin{equation} \label{eq:second_derivative}
2x'(s)y'(t) = \varphi ''(t - s),\;\; s,t \in T,  \; s<t,
\end{equation}
implying that $y'$ is zero on $T \cap (\tilde s,\infty )$ and thus that $x$ is strictly linearly increasing there, which is a contradiction to $x'(\tilde s)=0$. We conclude that $x'$ and by \eqref{eq:second_derivative} also $y'$ are nonzero on $T$. It now follows from \eqref{eq:second_derivative} that
\begin{equation} \nonumber
\frac{{x''(s)}}{{x'(s)}} =  - \frac{{y''(t)}}{{y'(t)}},\;\; s,t \in T,  \; s<t.
\end{equation}
Hence, $x'' - p x' = 0$ and $y'' + p y' =0$ on $T$, for some $p \in \mathbb{R}$. The rest of the proof is a straightforward verification. \hfill $\Box$

\begin{lemma} \label{lem3}
Suppose that $X$ is a nondeterministic $n$-parameter ($n \geq 1$) L\'evy process with exponent $\psi$ such that
$\psi (z) = a \psi (-z)$ for all $z \in \mathbb{R}^d$, for some $a\in\mathbb{R}$. Then $a=1$ and hence $X$ is symmetric.
\end{lemma}

\proof 
By the assumption, $ \psi (z) = a^2 \psi (z)$ for all $z \in \mathbb{R}^d$. The following facts complete the proof:
1) $X$ is deterministic if and only if $\psi (z) = - \psi ( - z)$ for all $z \in \mathbb{R}^d$ (consider the 
difference of two i.i.d. ID random variables); 2) $X$ is symmetric if and only if $ \psi (z) =  \psi (-z)$ for all $z \in \mathbb{R}^d$. \hfill $\Box$ \\

\begin{theorem} \label{thm2}
Let $X$ be a nondeterministic TPLP on $\mathbb{R}^d$ with exponent $\psi$, and $\alpha (t)=(x(t),y(t))_{t \in T}$ a decreasing path. If $X$ is symmetric, then $X^\alpha$ has stationary increments if and only if the path $\alpha$ meets one of conditions {\rm (i)-(iv)} in Lemma \ref{lem2}. For any of these conditions,
\begin{equation} \label{eq:thm_a}
{\rm E}[{\rm e}^{{\rm i} \langle z,X_t^\alpha  - X_s^\alpha  \rangle } ] =
\exp[{\varphi (t - s)\psi (z)}],\;z \in \mathbb{R}^d ,\;s,t \in T , \;s < t ,
\end{equation}
where $\varphi$ (here and below) is the corresponding function from Lemma \ref{lem2}. 
If on the other hand $X$ is not symmetric, then $X^\alpha$ has stationary increments if and only if $\alpha$ meets one of the conditions {\rm (i)} and {\rm (iv)}. For condition {\rm (iv)}, 
\begin{equation} \label{eq:thm_b}
{\rm E}[{\rm e}^{{\rm i} \langle z,X_t^\alpha  - X_s^\alpha  \rangle } ] = 
\exp \bigg[\varphi (t - s)\frac{{\psi (z) + \psi ( - z)}}{2}\bigg],\;s,t \in T , \;s < t .
\end{equation}
\end{theorem}

\proof
The first part of the theorem follows straight from \eqref{eq:difsym} and Lemma \ref{lem2}.
Assume therefore that $X$ is not symmetric. Then \eqref{eq:dif} implies that $X^\alpha$ has stationary increments if and only if 
\begin{equation} \label{eq:basic}
[m(B_{s,t}^{\rm l} ) - m(B_{s',t'}^{\rm l} )]\psi (z) = [m(B_{s',t'}^{\rm u} ) - m(B_{s,t}^{\rm u} )]\psi ( - z)
\end{equation}
for all $s,t,s',t' \in T$ such that $s<t$ and $t'-s'=t-s$, and all $z \in \mathbb{R}^d$. It follows straightforwardly from Lemma \ref{lem3} that \eqref{eq:basic} holds if and only if 
\begin{equation} \label{eq:two_equations}
m(B_{s,t}^{\rm l} ) = \eta _1 (t - s)\;\;{\rm and} \;\;m(B_{s,t}^{\rm u} ) = \eta _2 (t - s)
\end{equation}
for some functions $\eta_1$ and $\eta_2$; upon summation we see that $\alpha$ satisfies \eqref{eq:func_eq} with $\varphi=\eta_1 + \eta_2$. Among the solutions of \eqref{eq:func_eq}, only (i) and (iv) satisfy \eqref{eq:two_equations} (for suitable $\eta _1,\eta _2$), and \eqref{eq:thm_b} is verified by substitution into \eqref{eq:dif}. \hfill $\Box$ \\

We define the notion of an {\it increasing path} by letting $y$ in Definition \ref{dfn2} be nondecreasing (rather than nonincreasing). TPLPs along increasing paths are, roughly speaking, merely time changes of one-parameter L\'evy processes. The following remark says roughly that TPLPs do not have stationary increments along two-piece monotone paths.

\begin{remark} \label{rm2}
Let $T$ be the union of two intervals $T_1$ and $T_2$ having exactly one point in common, $a=\max T_1=\min T_2$. Suppose that $(x(t),y(t))_{t \in T_1}$ is an increasing path and $(x(t),y(t))_{t \in T_2}$ is a decreasing path, or vice versa, and that there exist $s \in T_1$ and $t \in T_2$ such that $y(a)\notin \{y(s),y(t)\}$. Assume that $X$ is a nonzero TPLP on $\mathbb{R}^d$. Then, it is easy to verify using Theorem \ref{thm2} that the process $\{X_{x(t),y(t)}:t \in T\}$ does not have stationary increments. 
\end{remark}

\section{The interesting cases}

In this section we discuss the classes of stationary increment processes corresponding to (ii)-(iv) in Lemma \ref{lem2}. 

\subsection{The independent increments case}

Let $X$ be a {\it symmetric} TPLP on $\mathbb{R}^d$ and $\alpha$ a decreasing `first vertical, then horizontal' path parameterized as in (ii) of Lemma \ref{lem2}. Thus, $X^\alpha$ has stationary independent increments. However, the notion of independent increments here is in the weak sense indicated in Sect. 2. Nevertheless, for any fixed time $t_0 \in T$ with $t_0 < \sup T$, the process $Y$ defined by
\begin{equation} \label{eq:useful}
Y_t  = X_{t + t_0}^\alpha   - X_{t_0 }^\alpha 
\end{equation}
is a L\'evy process in law on the time interval $[0,\sup T - t_0)$ (recall that a L\'evy process in law need not be right continuous with left limits), with characteristic exponent which is $ac$ times that of $X$. Indeed, $Y_0=0$, and the stationary independent increments property follows immediately from that of $X^\alpha$; the assertion on the exponent follows readily from \eqref{eq:thm_a}. 

Though very simple, equation \eqref{eq:useful} is instructive. Consider for example the process $W^\alpha$ where $W$ is a real-valued Brownian sheet. Based on \eqref{eq:useful} and elementary formulas for BM, one can straightforwardly find explicit formulas for, e.g., the probabilities
${\rm P}(\sup \nolimits_{t_0  \le t \le t_1 } (W_t^\alpha   - W_{t_0 }^\alpha  ) \le z)$, $z>0$, and ${\rm P}(W_t^\alpha   = W_{t_0 }^\alpha \;\,  {\rm for}\;{\rm some} \;\,  t_1\leq t \leq t_2  )$, for any $t_0,t_1,t_2 \in T$ with $t_0 < t_1 < t_2$.

\subsection{The case of uniformly scattered cancelling jumps}

From now until the end of the proof of Proposition \ref{prop3} below, $\alpha$ will denote the following decreasing path (corresponding to (iii) of Lemma \ref{lem2}):
\begin{equation} \nonumber
\alpha (t) = (bt,c - ct/l)_{t \in [0,l]},\;\; b,\,c,\,l>0.
\end{equation}
 The following simple lemma is of fundamental importance in what follows.

\begin{lemma} \label{lem4}
Let $(\xi_1,\xi_2)$ be uniformly distributed on the triangle $\Delta$ with vertices $(0,0)$, $(0,c)$, and $(bl,0)$. Let $\tau_1$ and $\tau_2$ be determined by $b\tau_1 = \xi_1$ and $c-c\tau_2/l = \xi_2$. Then, $(\tau_1,\tau_2)$ is distributed as $(U_{(1)},U_{(2)})$ where $U_{(1)}$ and $U_{(2)}$ are order statistics from a uniform distribution on $(0,l)$.
\end{lemma}

\proof 
The path $\alpha$ connects the vertices $(0,c)$ and $(bl,0)$. Thus, by construction, $0\leq \tau_1 \leq \tau_2 \leq l$. For any $u_1,u_2 \in (0,l)$ with $u_1<u_2$, we have
\begin{equation} \nonumber
\begin{split}
{\rm P}(\tau_1  \in {\rm d}u_1 ,\, \tau_2  \in {\rm d}u_2 ) & = {\rm P}(\xi _1  \in b\,{\rm d}u_1 ,\, \xi _2  \in c - (c/l){\rm d}u_2 ) \\
& = (2/l^2 ){\rm d}u_1 {\rm d}u_2,
\end{split}
\end{equation}
where we have used $[b(c/l)]/[(bl)c/2] = 2/l^2$. The joint density of $\tau_1$ and $\tau_2$ is thus that of $U_{(1)}$ and $U_{(2)}$ in the statement of the lemma. \hfill $\Box$ 


Given a measure $Q$ on $\mathbb{R}^d$, we denote by $\tilde Q$ the (`dual') measure defined for Borel subsets $B$ of $\mathbb{R}^d$ by $\tilde Q(B) = Q(-B)$. The following result accounts for the title of this subsection.

\begin{proposition} \label{prop3}
Let $Y$ be a one-parameter purely non-Gaussian L\'evy process on $\mathbb{R}^d$ with L\'evy measure $\nu$ satisfying
$\int_{\{ |x| \le 1\} } {|x|\nu ({\rm d}x)}  < \infty$ and zero drift. Let $(\tau_i)$ be an enumeration (possibly empty) of the jumping times of 
$Y$ in the time interval $(0,l)$, and $(J_i)$ the corresponding jumps. Further, let $(\tau'_i)$ be an independent sequence of i.i.d. uniform $(0,l)$ variables. Define the process $Y'$ on the time interval $[0,l]$ by 
\begin{equation} \label{eq:rearrangement}
Y'_t = \sum\limits_i {\mathbf{1}_{[0,t]} (\tau'_i )J_i } ,\;\;0 \leq t \leq l.
\end{equation}  
$Y' \stackrel{{\rm law}}{=} Y$ on $[0,l]$. Finally, define the process $Z$ on $[0,l]$ by
\begin{equation} \nonumber
Z_t = Y_t - Y'_t,\;\;0 \leq t \leq l,
\end{equation}
and let $X$ be a purely non-Gaussian TPLP on $\mathbb{R}^d$ with L\'evy measure $\hat \nu$ given by
\begin{equation} \label{eq:relation_measures}
\hat \nu  = (bc)^{ - 1} (\nu  + \tilde \nu )
\end{equation}
and zero drift (hence $X$ is symmetric). Then, $Z \stackrel{{\rm law}}{=} X^\alpha$.
\end{proposition}

It is useful to note the relation $\nu  = (bc/2) \hat \nu$ corresponding to \eqref{eq:relation_measures} in case $\nu$ is symmetric. In particular, there is no one-to-one correspondence between the measures $\nu$ and $\hat \nu$ satisfying \eqref{eq:relation_measures}.  

\begin{remark} \label{rm:generalize}
We do not find it important to consider the case where small jumps need to be compensated, i.e. the case
$\int_{\{ |x| \le 1\} } {|x|\nu ({\rm d}x)}  = \infty$. 
\end{remark}

{\noindent{\it Proof of Proposition \ref{prop3}.}} We first note that by the assumptions on $Y$, $Y'$ is well-defined and has the same law as $Y$.
Moreover, the condition on $\nu$ holds for $\hat \nu$ as well, and hence the drift of $X$ is defined.

The case of $Y$ being the zero process is trivial; assume first that $Y$ is a CPP (on $\mathbb{R}^d$) with rate $\lambda>0$ and jump distribution $F\,(\neq \delta_0)$. It follows from \eqref{eq:relation_measures} that $X$ is a two-parameter CPP with mean number of jumps in $[0,1]\times[0,1]$ equal to $2\lambda/(bc)$ and jump distribution $\hat F = {\textstyle{1 \over 2}} (F + \tilde F)$. Thus, in particular, the number of jumps of $X$ in the triangle $\Delta$ of Lemma \ref{lem4} is Poisson distributed with mean $[(bl)c/2][2\lambda /(bc)] = \lambda l$, just like that of $Y$ in the time interval $[0,l]$. Given that $X$ has $k\,(\geq 1)$ jumps in $\Delta$, the jump locations are independent and uniformly distributed in $\Delta$. Each jump results in a pair of $\hat F$-distributed cancelling jumps in the process $X^\alpha$. Then, since $\hat F$ is symmetric, Lemma \ref{lem4} leads us to conclude that the process $X^\alpha$ can be represented in law as a difference of two CPPs on $[0,l]$ with rate $\lambda$ and $\hat F$-distributed jumps, say $V$ and $V'$, such that $V'$ is obtained from $V$ by an independent rearrangement of the jumping times. It is easy to conclude that the same holds for the process $Z$ in the proposition. Thus $Z \stackrel{{\rm law}}{=} X^\alpha$. 

Assume now that $\nu(\mathbb{R}^d)=\infty$. Let $(Z_n)$ and $(X_n)$ be the sequences of processes obtained from $Z$ and $X$, respectively, by truncating the jumps smaller than $1/n$ in absolute value. For each $n$, it follows from the previous step that $Z_n \stackrel{{\rm law}}{=} X_n^\alpha$. Hence, the same is true for the limits $Z$ and $X^\alpha$. \hfill $\Box$ 

The following corollary and the paragraph that follows its proof provide an illuminating view of the Brownian bridge. 

\begin{corollary} \label{cor3}
Fix $l>0$. Let $(Y^n)$ be a sequence of one-parameter (nonzero) real-valued CPPs with rate $n$ and jump distribution $F$ with first moment $\mu'_1$ and second moment $\mu'_2$. For each $n$, let $Y'^n$ be the CPP on $[0,l]$ obtained from $Y^n$ by an independent rearrangement of the jumping times (as in \eqref{eq:rearrangement}), and define the process $Z^n$ on $[0,l]$ by
\begin{equation} \nonumber 
Z_t^n  = \frac{{Y_t^n - Y'^n_t }}{{\sqrt {2\mu '_2 n} }},\;\;0 \leq t \leq l.
\end{equation}
Then, $Z^n$ converges in FDDs to the standard Brownian bridge on $[0,l]$.
Moreover, $Z^n$ can be written as a difference of two limiting BMs on $[0,l]$ (in the sense of weak convergence) with variance parameter $1/2$ each, as follows:
\begin{equation} \nonumber
Z_t^n  = \frac{{Y_t^n  - n\mu '_1 t}}{{\sqrt {2\mu '_2 n} }} - \frac{{Y'^n_t  - n\mu '_1 t}}{{\sqrt {2\mu '_2 n} }}.
\end{equation}
\end{corollary}

\proof We first prove the second part of the assertion. Since weak convergence of L\'evy processes reduces to weak convergence of the marginal distributions at $t=1$ (see e.g. \cite[Corollary VII.3.6]{jacod}), we actually need to show that $(Y_1^n  - n\mu '_1 )/\sqrt {2\mu '_2 n}$ converges in distribution to the N$(0,1/2)$ law. This follows  from the central limit theorem upon replacing $Y_1^n$ by a sum of $n$ independent copies of $Y_1^1$ and using the equalities ${\rm E}[Y_{1}^1] = \mu '_1$, ${\rm Var}[Y_{1}^1 ] = \mu '_2$.

For the first part, define $V^n=Y^n/\sqrt {2\mu '_2 n}$ and $V'^n=Y'^n/\sqrt {2\mu '_2 n}$. By letting $V^n$, $V'^n$, and $Z^n$ play the role of $Y$, $Y'$, and $Z$ in Proposition \ref{prop3}, respectively, we get $Z^n \stackrel{{\rm law}}{=} (X^n)^\alpha$ where $\alpha (t) = (t,1-t/l)_{t \in [0,l]}$ and $X^n$ is a two-parameter CPP with L\'evy measure  $\hat \nu_n$ given by $\hat \nu_n = \nu_n + \tilde \nu_n$, $\nu_n$ being the L\'evy measure of $V^n$. In fact, $\nu_n + \tilde \nu_n$ is the L\'evy measure of the CPP $V^n-V^{\# n}$ where $V^{\# n}$ is an independent copy of $V^n$. Thus $X_{1,1}^n \stackrel{{\rm d}}{=} V^n_1-V^{\# n}_1$. By the convergence of $(Y_1^n  - n\mu '_1 )/\sqrt {2\mu '_2 n}$ mentioned above, $V^n_1-V^{\# n}_1$ and hence $X_{1,1}^n$ converge in distribution to the N$(0,1)$ law. Thus the FDDs of $X^n$ converge to those of a standard Brownian sheet $W$ and, in particular, those of $(X^n)^\alpha$ and hence of $Z^n$ to those of $W^\alpha$. Hence we are done since $W^\alpha$ has the same law as a standard Brownian bridge on $[0,l]$. \hfill $\Box$ 

\begin{figure}
	\centering
	\includegraphics[height=50mm]{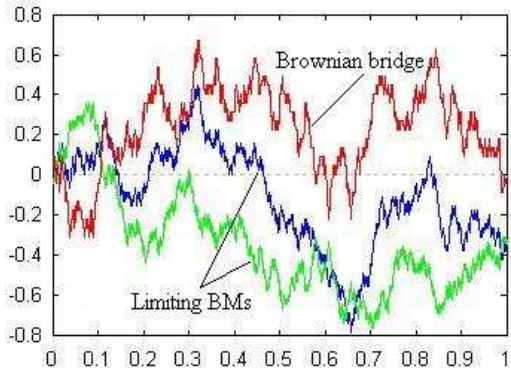}
	\caption {A standard Brownian bridge on $[0,1]$ viewed as a difference of two limiting BMs with variance parameter 1/2 each.}
	\label{fig2}
\end{figure}

Corollary \ref{cor3} invites us to consider a random walk analogy. It can be easily shown that if $\xi_1,\xi_2,\ldots$ are i.i.d. from a distribution $F$ as in the corollary and $(r_1,\ldots,r_{\left[ {nl} \right]})$ is a uniform random permutation of $(1,\ldots,{\left[ {nl} \right]})$, then the covariance function of the zero-mean process $\hat Z^n$ defined by
\begin{equation} \nonumber 
\hat Z_t^n  = \frac{{\sum\nolimits_{i = 1}^{\left[ {nt} \right]} {\xi _i }  - \sum\nolimits_{i = 1}^{\left[ {nt} \right]} {\xi _{r_i} } }}{{\sqrt {2\mu '_2 n} }},\;\;0 \leq t \leq l,
\end{equation}
is given by
\begin{equation} \nonumber
{\rm E}[\hat Z_s^n \hat Z_t^n ] =
\bigg(1 - \frac{{(\mu '_1 )^2 }}{{\mu '_2 }}\bigg)\frac{{\left[ {ns} \right]}}{n}\bigg(1 - \frac{{\left[ {nt} \right]}}{{\left[ {nl} \right]}}\bigg),\;\; 0 \leq s \leq t \leq l.
\end{equation}
Thus, the covariance function of $\hat Z^n$ converges as $n \rightarrow \infty$ to that of a Brownian bridge on $[0,l]$ with variance parameter 
$1 - (\mu '_1 )^2 /\mu '_2  \ge 0$. This motivates to consider convergence of $\hat Z^n$ to a Brownian bridge;
however, this problem lies well beyond our scope. The dynamics of $\hat Z^n$, $n=10^3$, for the basic example of simple symmetric random walk, 
${\rm P}(\xi _i  = \pm 1) = 1/2$, is illustrated in Figure \ref{fig2} (recall that in this case $\{ n^{ - 1/2} \sum\nolimits_{i = 1}^{\left[ {nt} \right]} {\xi _i } : t \in [0,1]\}$ converges weakly to a standard BM as $n \rightarrow \infty$); Monte Carlo simulations agreed well with pointwise convergence of $\hat Z^n_t$ to the N$(0,t(1-t))$ law. It is interesting to note \cite[p. 448]{biane} that Brownian bridge arises as the weak limit as $n \rightarrow \infty$ of a scaled simple symmetric random walk conditioned to come back to the origin after $n \in 2 \mathbb{N}$ steps ({\it random walk bridge}).

The Brownian `sheet-bridge relation' gives rise to the following remark.

\begin{remark} \label{rm:discrete}
Let $Y$ be a (nonzero) real-valued CPP with symmetric discrete jump distribution.
Let $\hat Y^{l}:=Y|Y_l=0$ be the associated `bridge process' on $[0,l]$. One can easily check/realize that $\hat Y^{l}$ has stationary increments. However, unlike in the brownian case, $\hat Y^{l}$ cannot be represented in law as $X^\alpha$, for any two-parameter CPP $X$ and decreasing path $\alpha$. By Theorem \ref{thm2}, we only need to consider the case $\alpha(t)=(bt,c-ct/l)_{t \in [0,l]}$. Assuming for a contradiction that $X^\alpha \stackrel{{\rm law}}{=} \hat Y^{l}$, it follows that $\sum\nolimits_{j = 0}^{n - 1} {{\mathbf 1}_{\{ X_{(j + 1)l/n}^\alpha - X_{jl/n}^\alpha   \ne 0\} } } \stackrel{{\rm d}}{=} \sum\nolimits_{j = 0}^{n - 1} {{\mathbf 1}_{\{ \hat Y_{(j + 1)l/n}^{l}  - \hat Y_{jl/n}^{l}  \ne 0\} } }$. Since the limits as $n \rightarrow \infty$, say $J^\alpha$ and $J^l$, are the number of jumps in $[0,l]$ of $X^\alpha$ and $\hat Y^{l}$, respectively, a contradiction will be reached if we show that ${\rm P}(J^\alpha   = 2k) \ne {\rm P}(J^l  = 2k)$ for some $k$. In fact, for some $p,q>0$ and all $k \in \mathbb{Z}_+$, ${\rm P}(J^\alpha   = 2k) = {\rm e}^{ - p } p ^k /k!$ but ${\rm P}(J^l  = 2k) \le q ^{2k} /(2k)!$. (In view of Remark \ref{rm:intuitive}, we note that $\hat Y^{l}$ is Markovian.)
\end{remark}

\subsection{The stationary case}

It follows readily from \eqref{eq:characteristic} that the stationary increment process $X^\alpha=\{ X_{a{\rm e}^{ct} ,b{\rm e}^{ - ct} } :t \in T\}$ $(a,b,c>0)$ is moreover (strictly) stationary, $X$ being an arbitrary TPLP on $\mathbb{R}^d$: for any $n \geq 1$ and points $t_j \in T$, $t_j+\tau \in T$, $j=1,\ldots,n$, we have $(X_{t_1 }^\alpha  , \ldots ,X_{t_n }^\alpha  ) \stackrel{{\rm d}}{=} (X_{t_1+\tau}^\alpha  , \ldots ,X_{t_n + \tau}^\alpha  )$. We can thus construct a stationary process whose one-dimensional marginal law is that of a given ID distribution on $\mathbb{R}^d$. 
If $X^\alpha$ is real-valued and square-integrable (${\rm E}[X_{1,1}^2 ] < \infty $), then its dependence structure is characterized as follows:
\begin{equation} \label{eq:correlation} 
\rho (u) = \frac{{{\rm Cov}(X_t^\alpha  ,X_{t + u}^\alpha  )}}{{\sqrt {{\rm Var}(X_t^\alpha  )} \sqrt {{\rm Var}(X_{t +  u}^\alpha  )} }} = {\rm e}^{-c|u|},\;\; t,\,t + u \in T.
\end{equation}
Indeed, for $u>0$ use the following decompositions into independent terms,
\begin{equation} \nonumber
\begin{split}
X_t^\alpha   & = X((0,a{\rm e}^{ct} ] \times (b{\rm e}^{ - c(t + u)} ,b{\rm e}^{ - ct} ]) +X ((0,a{\rm e}^{ct} ] \times (0,b{\rm e}^{ - c(t + u)} ]), \\
X_{t + u}^\alpha &  = X((0,a{\rm e}^{ct} ] \times (0,b{\rm e}^{ - c(t + u)} ]) + X ((a{\rm e}^{ct} ,a{\rm e}^{c(t + u)} ] \times (0,b{\rm e}^{ - c(t + u)} ]),
\end{split}
\end{equation}
and moreover the elementary fact that ${\rm Var}[X(B)] = m(B){\rm Var}[X_{1,1} ]$.

It is interesting to note that despite some resemblance, stationary processes of Ornstein--Uhlenbeck type do not belong to the class of stationary processes $X^\alpha$ if $X$ is non-Gaussian (i.e. having nonvanishing L\'evy measure). 
Let $V=\{V_t:t\geq 0\}$ be a stationary process of OU type on $\mathbb{R}^d$ (see \cite[Sect. 4]{barndorff3} and  \cite[Theorem 17.5]{sato}; see e.g. \cite{barndorff2,valdivieso} for the one-dimensional case), i.e. the stationary solution of a stochastic differential equation of the form ${\rm d}V_t  =  - cV_t\,{\rm d}t  + {\rm d}Z_{ct}$, 
where $Z$ is a L\'evy process with ${\rm E}[\log ^ +  |Z_1 |] < \infty$---termed the {\it background driving L\'evy process}---and $c>0$. The process $V$ is representable as 
\begin{equation} \label{eq:integral}
V_t = {\rm e}^{ - c t} V_0 + \int_0^t {{\rm e}^{ - c (t - s)} \,{\rm d}Z_{c s}}, 
\end{equation}
$V_0$ being independent of $Z$ and distributed as $\int_0^\infty  {{\rm e}^{ - cs}\,{\rm d}Z_{cs} }$.
The one-dimensional marginal law of the stationary process $V$ is an arbitrary {\it selfdecomposable} distribution, independent of $c$ owing to the unusual timing ${\rm d}Z_{ct}$. (We recall that selfdecomposable distributions constitute a very important class of ID distributions; see e.g. \cite{sato}.) Moreover, it is known \cite{barndorff2,valdivieso} and easy to check that if $V$ is real-valued and square-integrable, then its dependence structure is the same as in \eqref{eq:correlation}. However, we have the following result. 

\begin{proposition} \label{prop4}
Suppose that $V$ is a stationary process of OU type on $\mathbb{R}^d$, as in \eqref{eq:integral}, and that $X$ is a non-Gaussian TPLP such that $V_0 \stackrel{{\rm d}}{=} X_{1,1}$. Then the processes $V$ and $X^\alpha = \{ X_{{\rm e}^{c t} ,{\rm e}^{ - c t} } :t \ge 0\}$ are not identical in law.
\end{proposition}

In view of Remark \ref{rm:intuitive}, we note that processes of OU type (stationary or not) are Markovian (see \cite[Definition 17.2]{sato}).

{\noindent{\it Proof of Proposition \ref{prop4}.}}
Denote by $\psi$ the common log-characteristic function of $V_0$ and $X_{1,1}$, and by $\nu$ the corresponding L\'evy measure. Further, put $Q_t=\int_0^t {{\rm e}^{cs} \,{\rm d}Z_{cs} }$ and $R_t={\rm e}^{ct} X_t^\alpha  - X_0^\alpha$.
By \eqref{eq:integral}, ${\rm e}^{ct} V_t  = V_0  + \int_0^t {{\rm e}^{cs} \,{\rm d}Z_{cs} }$.
It thus suffices to show that if $Q_t \stackrel{{\rm d}}{=} R_t$ for all $t>0$, then $\nu(\mathbb{R}^d)=0$. From the equation for ${\rm e}^{ct} V_t$, by virtue of independence and stationarity, we find 
\begin{equation} \label{eq:LCF1}
{\rm E}[{\rm e}^{{\rm i}\langle z,Q_t \rangle } ] = \exp [\psi ({\rm e}^{ct} z) - \psi (z)],\;\; z \in \mathbb{R}^d,
\end{equation}
while upon decomposing $X_0^\alpha$ and $X_t^\alpha$ as we did below \eqref{eq:correlation}, we find 
\begin{equation} \label{eq:LCF2}
{\rm E}[{\rm e}^{{\rm i}\langle z,R_t \rangle } ] = \exp [{\rm e}^{ - ct} \psi (({\rm e}^{ct}  - 1)z) + (1 - {\rm e}^{ - ct} )(\psi ({\rm e}^{ct} z) + \psi ( - z))].
\end{equation}
Letting $t=t_n$ be such that ${\rm e}^{ct}  = n$, $n=2,3,\ldots\,$, and equating \eqref{eq:LCF1} and \eqref{eq:LCF2} (recall the first part of Remark \ref{rm:useful}), we obtain 
\begin{equation} \nonumber
\psi (nz) = \psi ((n - 1)z) + n\psi (z) + (n - 1)\psi ( - z).
\end{equation}
Then using induction we have
\begin{equation} \nonumber
\psi (nz) = 
\frac{{n(n + 1)}}{2}\psi (z) + \frac{{n(n - 1)}}{2}\psi ( - z).
\end{equation}
Interpreting the last equation in terms of (independent) L\'evy processes and in turn in terms of L\'evy measures, we conclude that
\begin{equation} \nonumber
\nu (n^{-1} B) = \frac{{n(n + 1)}}{2}\nu (B) + \frac{{n(n - 1)}}{2}\nu ( - B),
\end{equation}
for every Borel subset $B$ of $\mathbb{R}^d$. 
Letting $B_a  = \{ a/2 < |x| \le a\}$, we then get
\begin{equation} \nonumber
\int_{\{ 0 < |x| \le a\} } {|x|^2 \nu ({\rm d}x)}  \ge \sum\limits_{n = 0}^\infty  {\Big(\frac{a}{{2^{n + 1} }}\Big)^2  \nu (2^{-n} B_a )} \ge \sum\limits_{n = 0}^\infty  {\frac{{a^2 }}{8}\nu (B_a )}.
\end{equation}
Hence, by virtue of $\nu$ being a L\'evy measure, $\nu (B_a) = 0$ for all $a>0$. Thus clearly $\nu(\mathbb{R}^d)=0$, and so we are done. \hfill $\Box$ 

\subsection*{Acknowledgements}  

I am very happy to thank Mr and Mrs Shapack for funding my fellowship. 
I am indebted to my supervisor Prof. Ely Merzbach, whose guidance has inspired the present work.
I am grateful to two anonymous referees and an Associate Editor for their feedback and constructive criticism on an earlier draft of this paper.


\frenchspacing
\begin{thebibliography}{}

\bibitem{barndorff3}
Barndorff-Nielsen, O.~E., Maejima, M., and Sato, K.-I. (2006). Infinite divisibility for stochastic processes and time change. J. Theor. Probab. 19, 411--446.

\bibitem{barndorff2}
Barndorff-Nielsen, O.~E., and Shephard, N. (2001). Modelling by L\'evy processes for financial econometrics. In
L\'evy Processes: Theory and Applications (Barndorff-Nielsen, O.~E., Mikosch, T., Resnick, S., eds.),
Birkh\"auser, Boston, pp. 283--318.

\bibitem{biane} 
Biane, P., Pitman, J., and Yor, M. (2001). Probability laws related to the Jacobi theta and Riemann zeta functions, and Brownian excursions. Bull. Amer. Math. Soc. 38, 435--465.

\bibitem{dalang1}
Dalang, R.~C. (2003). Level sets and excursions of the Brownian sheet. In CIME 2001 summer school, Topics in Spatial Stochastic Processes (Merzbach, E., ed.), Lect. Notes Math. 1802, Springer, pp. 167--208.

\bibitem{dalang2}
Dalang, R.~C., and Walsh, J.~B. (1992). The sharp Markov property of L\'evy sheets. Ann. Probab. 20, 591--626.

\bibitem{di-nardo}
Di Nardo, E., Nobile, A.~G., Pirozzi, E., and Ricciardi, L.~M. (2001). A computational approach to first-passage-time problems for Gauss--Markov processes. Adv. Appl. Probab. 33, 453--482.

\bibitem{ehm}
Ehm, W. (1981). Sample function properties of multi-parameter stable processes. Z. Wahrsch. Verw. Gebiete 56, 195--228. 

\bibitem{jacod}
Jacod, J., and Shiryaev, A.~N. (2003). Limit Theorems for Stochastic Processes (2nd ed.), Springer, Berlin.

\bibitem{lagaize}
Lagaize, S. (2001). H\"older exponent for a two-parameter L\'evy process. J. Multivariate Anal. 77, 270--285.

\bibitem{sato}
Sato, K-I. (1999). L\'evy Processes and Infinitely Divisible Distributions, Cambridge University  Press. 

\bibitem{valdivieso} 
Valdivieso, L., Schoutens, W., and  Tuerlinckx, F. (2009). Maximum likelihood estimation in processes of Ornstein--Uhlenbeck type. Stat. Infer. Stoch. Process. 12, 1--19.

\bibitem{vares}
Vares, M.~E. (1982). Representation of the square integrable martingales generated by a two-parameter L\'evy process. Z. Wahrsch. Verw. Gebiete 61, 161--188.  


\end{thebibliography}
\end{document}